\documentclass[11pt,fleqn,twoside]{article}
\usepackage{amsfonts,latexsym}
\makeatletter
\newcommand{\prava}{\footnotesize\it
\begin{flushright}
\begin{minipage}{18cm}
Copyright \copyright 1998 by D.J. Arrigo and J.M. Hill
\end{minipage}
\end{flushright}}

\newcommand{\name}[1]{\begin{flushleft}
                       \LARGE \bf #1
                       \end{flushleft}\vspace{-3mm}}

\newcommand{\Author}[1]{\begin{flushleft}
                       \it #1 \end{flushleft}}

\newcommand{\Adress}[1]{\begin{flushleft}
                       \it #1 \end{flushleft}}

\newcommand{\Date}[1]{\begin{flushleft}
                      \small  \it #1 \end{flushleft}}

\newcommand{\ehkol}{Author \ name}
\newcommand{\ohkol}{Article \ name}
\renewcommand{\@evenhead}{
\hspace*{-3pt}\raisebox{-15pt}[\headheight][0pt]{\vbox{\hbox to \textwidth
{\thepage \hfil \ehkol}\vskip4pt \hrule}}}
\renewcommand{\@oddhead}{
\hspace*{-3pt}\raisebox{-15pt}[\headheight][0pt]{\vbox{\hbox to \textwidth
{\ohkol \hfil \thepage}\vskip4pt\hrule}}}
\renewcommand{\@evenfoot}{}
\renewcommand{\@oddfoot}{}

\newcommand{\be}{\begin{equation}}
\newcommand{\ee}{\end{equation}}
\newcommand{\ba}{\hspace*{-5pt}\begin{array}}
\newcommand{\ea}{\end{array}}

\makeatother

\begin{document}
\setcounter{page}{115}

\renewcommand{\theequation}{\arabic{section}.\arabic{equation}}
\renewcommand{\thesection}{\arabic{section}}

\thispagestyle{empty}

\renewcommand{\ehkol}{D.J. Arrigo and  J.M. Hill}
\renewcommand{\ohkol}{On a Class of Linearizable Monge-Amp\`ere Equations}

\begin{flushleft}
\footnotesize {\sf
Journal of Nonlinear Mathematical Physics \qquad 1998, V.5, N~2},\ 
\pageref{arrigo-fp}--\pageref{arrigo-lp}.
\hfill
{{\sc Letter}}
\end{flushleft}

\vspace{-5mm}

\renewcommand{\footnoterule}{}
{\renewcommand{\thefootnote}{}
 \footnote{\prava}}

\name{On a Class of Linearizable \\ Monge-Amp\`ere Equations}
\label{arrigo-fp}

\Author{D.J. ARRIGO and J.M. HILL}

\Adress{Department of Mathematics, University of Wollongong,\\
Wollongong, NSW 2522, Australia}

\Date{Received November 21, 1997; Accepted November 25, 1997}

\begin{abstract}
\noindent
Monge-Amp\`ere equations of the form,
$u_{xx}u_{yy}-u_{xy}^2=F(u,u_x,u_y)$ arise in many areas of f\/luid and
solid mechanics. Here it is shown that in the special case
$F=u_y^4f(u, u_x/u_y)$, where $f$ denotes an arbitrary function, the
Monge-Amp\`ere equation can be linearized by using a sequence of
Amp\`ere, point, Legendre and rotation transformations. This
linearization is a generalization of three examples from f\/inite
elasticity, involving plane strain and plane stress deformations of
the incompressible perfectly elastic Varga material and also relates
to a previous linearization of this equation due to Khabirov~[7].
\end{abstract}

\section{Introduction}

Monge-Amp\`ere equations of the form
\be
u_{xx}u_{yy} -u_{xy}^2=F(x,y,u,u_x,u_y),
\ee
are well known to arise is many areas of science and engineering,
especially areas relating to f\/luid mechanics (see for example von
Mises [1] and Martin [2]). In three recent papers Hill and Arrigo [3,
4] and Arrigo and Hill [5], the present authors have shown that
Monge-Amp\`ere equations also arise in the context of f\/inite elastic
deformations. In particular, it is shown that certain plane strain,
plane stress and axially symmetric deformations of the incompressible
perfectly elastic Varga material (see Varga [6]) all give rise to
Monge-Amp\`ere equations
of the form (1.1). Three of these may by linearized by a
sequence of elementary transformations and the purpose of this brief
communication is to show that the same sequence of transformations is
also ef\/fective when $F=u^4_yf(u,u_x/u_y)$, where $f$ denotes an
arbitrary function.

The recent work of the authors in f\/inite elasticity may summarized as
follows. For plane strain deformations of the form
\be
x=x(X,Y), \qquad y=y(X,Y),
\ee
where $(X,Y)$ and $(x,y)$ denote respectively material and spatial
plane rectangular Cartesian coordinates Hill and Arrigo [3, 4] show
that for the Varga elastic material, exact f\/inite elastic
deformations can be determined from the following
\be
x=U_X, \qquad y=U_Y,
\ee
\be
x=V_\alpha, \qquad y=V_\beta,
\ee
\be
x=\frac{W_X}{W_X^2+W^2_Y}, \qquad y=-\frac{W_Y}{W^2_X+W^2_Y},
\ee
where $U(X,Y)$, $V(\alpha,\beta)$ and $W(X,Y)$ all satisfy Monge-Amp\`ere
equations, namely
\be
U_{XX}U_{YY}-U_{XY}^2=1,
\ee
\be
V_{\alpha\alpha}V_{\beta\beta}-V_{\alpha\beta}^2=(\alpha^2+\beta^2)^{-2},
\ee
\be
W_{XX}W_{YY}-W_{XY}^2=(W_X^2+W_Y^2)^2,
\ee
where $(\alpha, \beta)$ are intermediate coordinates which are
def\/ined by
\be
\alpha=\frac{X}{X^2+Y^2}, \qquad \beta=-\frac{Y}{X^2+Y^2}.
\ee
In addition, for axially symmetric deformations of the form
\be
r=r(R,Z), \qquad \theta =\Theta, \qquad z=z(R,Z),
\ee
where $(R,\Theta, Z)$ and $(r,\theta,z)$ denote respectively material
and spatial cylindrical polar coordinates, Hill and Arrigo [3] and Arrigo and
Hill [5] show that for the Varga elastic material, exact f\/inite
elastic deformations may be determined from
\be
r=U_R, \qquad z=U_Z,
\ee
\be
r=V_\alpha, \qquad z=V_\beta,
\ee
where $U(R,Z)$ and $V(\alpha,\beta)$ satisfy respectively the
Monge-Amp\`ere equations
\be
U_{RR}U_{ZZ}-U_{RZ}^2=\frac{R}{U_R},
\ee
\be
V_{\alpha\alpha}V_{\beta\beta}-V_{\alpha\beta}^2=\frac{\alpha}
{(\alpha^2+\beta^2)^2V_\alpha},
\ee
where $(\alpha, \beta)$ are now the intermediate coordinates def\/ined
by
\be
\alpha=\frac{R}{R^2+Z^2}, \qquad \beta=-\frac{Z}{R^2+Z^2}.
\ee

Finally, for plane stress deformations arising from the membrane
theory of thin plane elastic sheets, which assumes a three
dimensional deformation of the form
\be
x=x(X,Y), \qquad y=y(X,Y), \qquad z=\lambda(X,Y)Z,
\ee
Arrigo and Hill [5] show that f\/inite elastic solutions may be
determined from
\be
x=V_\alpha, \qquad y=V_\beta,
\ee
where $\alpha$ and $\beta$ are precisely as def\/ined in (1.9) and
$V(\alpha, \beta)$ in this case satisf\/ies the Monge-Amp\`ere equation
\be
V_{\alpha\alpha}V_{\beta\beta}-V^2_{\alpha\beta}=
(\alpha^2+\beta^2)^{-2}(\alpha V_\alpha +\beta V_\beta-V)^{-1}.
\ee
Under the Legendre transformation
\[
u=\alpha V_\alpha +\beta V_\beta-V, \qquad x=V_\alpha, \qquad
y=V_\beta,
\]
the above Monge-Amp\`ere equations (1.7) and (1.18) can all be
embodied in the general equation
\be
u_{xx}u_{yy}-u_{xy}^2=A(u_x^2+u_y^2)^2,
\ee
where $A=1$ for equation (1.7) and $A=u$ for equation (1.18) and
of course $A=1$ for equation (1.8). In Hill and Arrigo [4] and Arrigo
and Hill [5] these special cases of equation (1.19), namely (1.7),
(1.8) and (1.18) are shown to be linerazable under a combination of a
Amp\`ere, point, Legendre and rotational transformations. Remarkably,
these Monge-Amp\`ere equations belong to a much larger class of
Monge-Amp\`ere equations which can be linearized by precisely the
same combination of transformations. The purpose of this paper is to
show that the Monge-Amp\`ere equation
\be
u_{xx}u_{yy}-u^2_{xy}=u^4_yf(u,u_x/u_y),
\ee
where $f$ denotes an arbitrary function, may be transformed to the
linear equation
\be
U_{XX}+f(X,Y)U_{YY}=0,
\ee
under the contact transformation (2.2) and this result is established
in the following section.

\setcounter{equation}{0}

\section{The basic linearization}

Following Hill and Arrigo [4] and Arrigo and Hill [5], we consider
the following sequence of transformations applied to the
Monge-Amp\`ere equation (1.1), namely
\be
\ba{ll}
\mbox{(i)} & x=\alpha, \quad y=V_\beta, \quad u=V-\beta V_\beta, \\
\mbox{(ii)} & \alpha=\xi, \quad \beta= 1/\eta, \quad V=W/\eta,\\
\mbox{(iii)} & \xi=Z_\tau, \quad \eta=Z_\sigma, \quad W=\tau Z_\tau+\sigma
Z_\sigma -Z,\\
\mbox{(iv)} & \tau=-Y, \quad \sigma =X, \quad Z=-U.
\ea
\ee
The f\/irst transformation represents an Amp\`ere transformation, the
second is a simple point transformation, the third is a Legendre
transformation, while the fourth represents a rotation and scaling
transformation. Combining all four transformations we have
\be
x=U_Y, \qquad y=U-YU_Y, \qquad u=X.
\ee
Now from (2.2), we f\/ind the f\/irst and second order partial
derivatives according to the following
\be
\ba{l}
u_x=Y/U_X, \qquad u_y=1/U_X,\\[1mm]
u_{xx}=\left(Y^2U^2_{XY}-Y^2U_{XX}U_{YY}-2YU_XU_{XY}+U^2_{X}\right)
U_X^{-3}U_{YY}^{-1},\\[1mm]
u_{xy}=\left(YU_{XY}^2-YU_{XX}U_{YY}-U_XU_{XY}\right)
U_X^{-3}U_{YY}^{-1},\\[1mm]
u_{yy}=\left(U_{XY}^2-U_{XX}U_{YY}\right)U^{-3}_XU_{YY}^{-1},
\ea
\ee
which, on substitution into equation (1.1), and assuming that
the function $F$ is independent of $x$ and $y$, yields the equation
\be
U_{XX}+F(X,Y/U_X, 1/U_X)U_X^4U_{YY}=0.
\ee
If now we require that equation (2.4) be linear, namely
\be
U_{XX}+f(X,Y)U_{YY}=0,
\ee
for some function $f(X,Y)$, then we require that the following equation
\be
F(X,Y/U_X,1/U_X)U_X^4=f(X,Y),
\ee
holds identically. Now on transforming back to the original variables
$(x,y,u)$, using (2.1) and (2.3) we may deduce
\be
F(u,u_x,u_y)=u_y^4f(u,u_x/u_y),
\ee
which is the desired result.

We note that under only the Amp\`ere transformation (i) of (2.1),
that the Monge-Amp\`ere equation (1.1) with $F=1$ is linearizable to
Laplaces' equation. We also comment that Khabirov [7] showed, using
Lie contact symmetry analysis, that Monge-Amp\`ere equations of the
form
\be
u_{xx}u_{yy}-u_{xy}^2=F(x,y),
\ee
can be linearized provided that $F(x,y)$ is one of the following four
cases, namely
\[
F=0, \quad F=1, \quad F=f(x), \quad F=x^{-4}g(y/x),
\]
where $f(x)$ and $g(y/x)$ denote arbitrary functions. From our point
of view the f\/inal case is of particular interest since under the
Legendre transformation,
\be
x=U_X, \quad y=U_Y, \quad u=XU_X+YU_Y-U,
\ee
the Monge-Amp\`ere equation
\be
u_{xx}u_{yy}-u_{xy}^2=x^{-4}g(y/x)=y^{-4}g^*(y/x),
\ee
transforms to the equation
\be
U_{XX}U_{YY}-U_{XY}^2=U_Y^4G^*(U_Y/U_X),
\ee
where $G^*=1/g^*$ and this equation is of the form (1.1) where $F$ is
as given in (2.7) except that in this case $f$ is independent of $u$.
Although this is a special case of the results presented here,
Khabirov [7] does not provide the explicit contact transformation
which linearizes equation (2.11).

\label{arrigo-lp}

\end{document}